\documentclass[preprint,12pt]{elsarticle}

\usepackage{amssymb}
\usepackage{amsmath}
\usepackage{amsthm}

\usepackage{mathtools}

\newcommand{\N}{\mathbb{N}}
\newcommand{\R}{\mathbb{R}}
\newcommand{\C}{\mathbb{C}}
\newcommand{\bz}{\mathbf{z}}
\newcommand{\bw}{\mathbf{w}}
\newcommand{\by}{\mathbf{y}}
\newcommand{\bx}{\mathbf{x}}

\DeclareMathOperator{\Id}{Id}
\DeclareMathOperator{\Tr}{Tr}
\DeclareMathOperator{\Ker}{Ker}
\DeclareMathOperator{\Imm}{Im}
\DeclareMathOperator{\Ree}{Re}
\DeclareMathOperator{\LEX}{LEX}
\DeclareMathOperator{\LT}{LT}
\DeclareMathOperator{\ch}{ch}

\theoremstyle{definition}
\newtheorem{Def}{Definition}[section]
\newtheorem{ex}[Def]{Example}
\newtheorem{rmk}[Def]{Remark}

\theoremstyle{plain}
\newtheorem{prop}[Def]{Proposition}
\newtheorem{coro}[Def]{Corollary}
\newtheorem{teo}[Def]{Theorem}
\newtheorem{lemma}[Def]{Lemma}
\newtheorem{con}[Def]{Conjecture}

\journal{Linear Algebra and its Applications}

\begin{document}

\begin{frontmatter}



\title{The Hermitian Killing form and root counting of complex polynomials with conjugate variables}


\author{Davide Furchì} %% Author name

\affiliation{organization={Dipartimento di
		Scienza e Alta Tecnologia, Università dell'Insubria},%Department and Organization
            %addressline={}, 
            city={Como},
            postcode={22100}, 
            state={},
            country={Italy}}

%% Abstract
\begin{abstract}
Inspired by the work about solutions of a system of real polynomial equations done by Hermite, this paper introduces a Hermitian form, which encodes information about solutions of a system of complex polynomial equations with conjugate variables.

Adopting the presented object, a new general bound for the number of solutions of a harmonic polynomial equation is proved.
\end{abstract}


%\item A different method to count the number of solutions of a system of polynomial equations with conjugate variables
%\item New upper bound for the number of zeros of Harmonic polynomial equations

%% Keywords

\end{frontmatter}



\section{Introduction}\label{sec:introduction}

In many theoretical and practical contexts it happens to perform calculation with the conjugate operation. In particular, we are interested in the zeros of \emph{generalized polynomials}, i.e.\ complex polynomials, where the variables could be conjugated.

Let $\bz=(z_1,\ldots,z_r)$, if we restrict our attention to generalized polynomials of the form
\begin{equation*}
  h(\bz,\bar{\bz})=q(\bz)+p(\bar{\bz}),
\end{equation*}
where $q$ and $p$ are complex multivariate polynomials in $\bz$ and $\bar{\bz}$ respectively, we get the special case of harmonic polynomial equations, which are of particular interest, e.g.\ see \cite{kls2018,ll2016,bbdhs2020,hllm2014}.

Let $r=1$ and consider
\begin{equation}\label{harmpol}
 q(z)=b_mz^m+\ldots+b_1z\qquad\text{and}\qquad p(\bar{z})=\bar{z}^n+a_{n-1}\bar{z}^{n-1}+\ldots+a_1\bar{z}+a_0 
\end{equation}
complex univariate generalized polynomials where the degree zero term is given by the coefficient $a_0$ of $h$. The B\'ezout Theorem gives the naive bound $n^2$ for the number of isolated roots of a harmonic polynomial $h(z,\bar{z}
)=q(z)+p(\bar{z})$ with $m=\deg(q)\leq\deg(p)=n$. In \cite{w1998} this bound was proved to be sharp for $m=n$, and the general bound
\begin{equation}\label{wil}
    3n-2+m(m-1)
\end{equation}
was conjectured for $m<n$.

For $m=n-1$, formula \eqref{wil} yields the bound $n^2$, whose sharpness was proved in \cite{w1998,bhs1995}.

For $m=1$, we get the bound $3n-2$, which was shown to be correct for any $n$ and even sharp for $n=2,3$ in \cite{ks2002}. Then, this bound was shown to be sharp also for $n=4,5,6,7,8$ in \cite{bl2004} and for any $n$ in \cite{l2008}.

However, in \cite{lll2013} it was shown that the bound \eqref{wil} is in general false by providing counterexamples for $m=n-3$ and was proposed the general bound
\begin{equation}\label{lll}
    2m(n-1)+n.
\end{equation}

Formula \eqref{lll} is linear in $n$ and is equal to the sharp bound $3n-2$ for $m=1$. However, if $m\geq n/2$, formula \eqref{lll} is greater than the naive bound $n^2$. One result of this work is the following, which sharpens the current bound $n^2$ for $1<m<n-1$.

\begin{teo}\label{main}
Using the notations of equations \eqref{harmpol}, if $n-2\geq m$ the harmonic equation $h(z,\bar{z})=q(z)+p(\bar{z})=0$ admits at most $n^2-1$ solutions if $(n-1)a_{n-1}^2-2na_{n-2}=0$ and at most $n^2-2$ solutions otherwise.
\end{teo}

In the general setting of equations of generalized polynomials, the present literature suggests the use of Gr\"obner basis as in \cite{k1987}, or to subdivide the polynomials into their real and imaginary part and study the problem in a real setting, for example by using the Killing form. A first work on the Killing form is \cite{prs1993}, consult \cite{pv2020} for a comprehensive dissertation and \cite{ss2002} for an introduction on harmonic polynomials.

We define a new method that permits to avoid this procedures and yields a more manageable object. Other than the theoretical interest of this theory, the method proposed is of computational importance. Admitting conjugation, the Hermitian Killing form requires half the variables and the coefficient parameters, compared to the Killing form applied to the associated real system.

%------------------------------------------------------------------------------------------

\section{Preliminaries and notations}

We will work in general in a complex setting. An ideal generated by polynomials is denoted with $\langle\ldots\rangle$, the relative set of common zeros will be indicate with $V(\ldots)$ and the set of leading terms of an ideal will be denoted with $\LT(\cdot)$. Let $r\in\N$, for $\bx\in\C^r$ and $\alpha\in\N^{r}$ we denote with $\bx^{\alpha}$ the monomial $x_1^{\alpha_1}\cdots x_r^{\alpha_r}$. The symbol $|\cdot |$ indicates the classical module in $\C$. Throughout the work, we assume we have set a monomial order. For computations we use the $\LEX$ monomial order. The superscripts $T$ and $H$ used for matrices will denote the transpose and the conjugated transpose, respectively. The symbol $i$ will stand for the imaginary unit.

\vspace{1cm}
Let $\bx=(x_1,\ldots,x_r)$ and consider a zero dimensional ideal
\begin{equation*}
I=\langle p_1(\bx),\ldots,p_c(\bx)\rangle\subset\C[\bx]
\end{equation*}
for $c\in\N$. In particular, we have $V(I)=\lbrace\bx^{(1)},\ldots,\bx^{(d)}\rbrace\subset\C^r$ for some $d>0$. The multiplications by $x_k$ for $k=1,\ldots,r$ induce linear applications on the quotient $\C[\bx]/I$ as follows. For $[f]\in\C[\bx]/I$ we define
\begin{equation*}
  M_{x_k}([f])\coloneqq [fx_k].
\end{equation*}
The matrices representing these maps, with respect to the basis $\lbrace [\bx^{\alpha}]\mid\bx^{\alpha}\notin\LT(I)\rbrace$, are the \emph{companion matrices} of $I$. More generally, for $[g]\in\C[\bx]/I$ we define the product $M_{g}([f])\coloneqq [fg]$.

We have the following results, see \cite[Chapter 2, Section 4]{clo2002}.

\begin{lemma}\label{prodeigen}
Let $f,g\in\C[\bx]$, then 
\begin{enumerate}[\normalfont i)]
    \item the map $M_f$ is well defined,
    \item it holds $M_{f}+M_{g}=M_{f+g}$,
    \item it holds $M_{f}M_{g}=M_{g}M_{f}=M_{fg}$,
    \item it holds $M_{f}=f(M_{x_1},\ldots,M_{x_n})$,
    \item the eigenvalues of $M_{f}$ are the evaluations $f(\bx^{(k)})\in\C$ of $\bx^{(k)}\in V(I)$ for $k=1,\ldots,d$ with relative multiplicities.
\end{enumerate} 
\end{lemma}

Using this definition we have the following results, see \cite[Chapter 4, Section 2]{clo2002}.

\begin{prop}\label{decA}
Let $\xi\in\C[\bx]$ be such that the evaluations $\xi(\bx^{(k)})$ are all different for $k=1,\dots,d$.  Consider for any $k=1,\ldots,d$ the linear applications $M_{\xi(\bx)-\xi(\bx^{(k)})}$. Then, setting $A_k=\bigcup_{n\in\N}\Ker(M_{\xi(\bx)-\xi(\bx^{(k)})})^{n}$, we have the decomposition
\begin{equation*}
    \C[\bx]/I=\oplus_{k=1}^dA_k.
\end{equation*}
If $f\in A_k$ then $f(\bx^{(j)})\neq 0$ for $k\neq j$. Any sub algebras $A_k$ admits a unity $e_k$ such that $e_k^2=e_k$, $e_ke_j=0$ for $k\neq j$ and $e_k(\bx^{(j)})=\delta_{k,j}$ for any $k,j$. Moreover, it holds $\bar{e}_k=e_j$ if $\bar{\bx}^{(k)}=\bx^{(j)}$.
\end{prop}

%-------------------------------------------------------------------------------------------

\section{From the real to the complex case}\label{herm}

\begin{Def}
A \emph{generalized polynomial} is an expression of the form 
\begin{equation*}
    p(\bz,\bar{\bz})=\sum_{(\alpha,\beta)}a_{\alpha,\beta}\bz^{\alpha}\bar{\bz}^{\beta},
\end{equation*}
where the coefficients $a_{\alpha,\beta}\in\C$ are non zero only for finitely many multi-indices $(\alpha,\beta)\in\N^{2r}$. We denote its conjugate $\bar{p}(\bz,\bar{\bz})=\sum_{(\alpha,\beta)}\bar{a}_{\alpha,\beta}\bz^{\alpha}\bar{\bz}^{\beta}$.
\end{Def}

We are interested in solving equations of the type $p(\bz,\bar{\bz})=0$.

%-------------------------------------------------------------------------------------------

\subsection{The Killing or Trace form}
Let us recall the construction of the Killing form. Consider a zero dimensional ideal
\begin{equation*}
I=\langle p_1(\bx),\ldots,p_c(\bx)\rangle\subset\R[\bx]
\end{equation*}
for $c\in\N$. In particular, we have
\begin{equation*}
V(I)=\lbrace\bx^{(1)}\ldots,\bx^{(d)}\rbrace\subset\C^r
\end{equation*}
for some $d>0$. Denote with $m_k$ the multiplicity of $\bx^{(k)}$ for $k=1,\ldots,d$, then it holds $\dim(\C[\bx]/I)=\sum_{k=1}^dm_k$.

\begin{Def}
We define for $\xi\in\R[\bx]$ the symmetric form
\begin{equation*}
  \C[\bx]/I\times\C[\bx]/I\to\C,\qquad([f],[g])\mapsto\mathcal{K}_{\R}^{\xi}(f,g)\coloneqq\Tr(M_{\xi}M_fM_{g}).
\end{equation*}
The form $\mathcal{K}_{\R}^{1}$ is called the \emph{Killing} or \emph{Trace} form.  
\end{Def}

The following result is due to Hermite, see \cite[Chapter 2, Theorem 5.2]{clo2002}.

\begin{teo}\label{hermite}
Let $\xi\in\R[\bx]$ and consider the restriction of $\mathcal{K}^\xi_{\R}$ on $\R[\bx]/I\times\R[\bx]/I$, then:
\begin{enumerate}[\normalfont i)]
    \item The variety $V(I)$ consists only in real points $\bx^{(k)}$ with multiplicity $1$ such that $\xi(\bx^{(k)})>0$ iff $\mathcal{K}_{\R}^{\xi}$ is positive definite.
    \item The rank of $\mathcal{K}_{\R}^{\xi}$ is the number of distinct points $\bx^{(k)}\in V(I)$ for which $\xi(\bx^{(k)})\neq 0$.
    \item The number of distinct real points $\bx^{(k)}\in V(I)$ such that $\xi(\bx^{(k)})>0$ minus the number of distinct real points $\bx^{(k)}\in V(I)$ such that $\xi(\bx^{(k)})<0$ equals the signature of $\mathcal{K}_{\R}^{\xi}$.
\end{enumerate}
If moreover $\xi$ is such that the evaluations $\xi(\bx^{(k)})$ are all different for $k=1,\ldots,d$, then:
\begin{enumerate}[\normalfont i)]
    \setcounter{enumi}{3}
    \item The number of distinct real points $\bx^{(k)}\in V(I)$ such that $\xi(\bx^{(k)})>0$ equals the number of positive eigenvalues of $\mathcal{K}_{\R}^{\xi}$ minus the number of negative eigenvalues of $\mathcal{K}_{\R}^1$.
    \item The number of distinct real points $\bx^{(k)}\in V(I)$ such that $\xi(\bx^{(k)})<0$ equals the number of negative eigenvalues of $\mathcal{K}_{\R}^{\xi}$ minus the number of negative eigenvalues of $\mathcal{K}_{\R}^1$.
\end{enumerate}
\end{teo}

%-------------------------------------------------------------------------------------------

\subsection{The Hermitian Killing form}
Set $c\in\N$ and note that when solving a system of generalized polynomial equations $p_k(\bz,\bar{\bz})=0$ for $k=1,\ldots,c$, we also can state $\bar{p}_k(\bar{\bz},\bz)=0$ for $k=1,\ldots,c$. Thus, we introduce the variables $\bw=(w_1,\ldots,w_r)$, where $w_k$ represents $\bar{z}_k$ for $k=1,\ldots,r$ and consider a zero dimensional ideal 
\begin{equation}\label{conjideal}
\tilde{I}=\langle p_1(\bz,\bw),\ldots,p_c(\bz,\bw),\bar{p}_1(\bw,\bz),\ldots,\bar{p}_c(\bw,\bz)\rangle\subset\C[\bz,\bw].
\end{equation}
In particular, we have 
\begin{equation*}
    V(\tilde{I})=\lbrace (\bz^{(1)},\bw^{(1)}),\ldots,(\bz^{(d)},\bw^{(d)})\rbrace\subset\C^{2r}
\end{equation*}
for some $d>0$. Denote with $m_k$ the multiplicity of $(\bz^{(k)},\bw^{(k)})$ for $k=1,\ldots,d$, then it holds $\dim(\C[\bz,\bw]/\tilde{I})=\sum_{k=1}^dm_k$.

\begin{lemma}\label{associated}
Consider a polynomial $p(\bz,\bw)\in\C[\bz,\bw]$, a point $(\bz^{(\lambda)},\bw^{(\lambda)})\in\C^{2r}$ and $\lambda\in\C$, then it holds $p(\bz^{(\lambda)},\bw^{(\lambda)})=\bar{p}(\bw^{(\lambda)},\bz^{(\lambda)})=\lambda$ iff $p(\bar{\bw}^{(\lambda)},\bar{\bz}^{(\lambda)})=\bar{p}(\bar{\bz}^{(\lambda)},\bar{\bw}^{(\lambda)})=\bar{\lambda}$. Moreover, if $\bz^{(\lambda)}=\bar{\bw}^{(\lambda)}$ and $\lambda\in\R$ then any of the former equalities implies the other three.
\end{lemma}
\begin{proof}
By conjugation, $p(\bz^{(\lambda)},\bw^{(\lambda)})=\bar{p}(\bw^{(\lambda)},\bz^{(\lambda)})=\lambda$ iff 
\[
p(\bar{\bw}^{(\lambda)},\bar{\bz}^{(\lambda)})=\overline{\bar{p}(\bw^{(\lambda)},\bz^{(\lambda)})}=\bar{\lambda}\qquad\text{and}\qquad\bar{p}(\bar{\bz}^{(\lambda)},\bar{\bw}^{(\lambda)})=\overline{p(\bz^{(\lambda)},\bw^{(\lambda)})}=\bar{\lambda}.
\]
If $\bz^{(\lambda)}=\bar{\bw}^{(\lambda)}$ then 
\[
p(\bar{\bw}^{(\lambda)},\bar{\bz}^{(\lambda)})=p(\bz^{(\lambda)},\bar{\bz}^{(\lambda)})=\overline{\bar{p}(\bar{\bz}^{(\lambda)},\bz^{(\lambda)}})=\overline{\bar{p}(\bw^{(\lambda)},\bz^{(\lambda)}})=\lambda
\]
and similarly for the other two.
\end{proof}

\begin{lemma}\label{asso}
If $(\bz^{(k)},\bw^{(k)})\in V(\tilde{I})$ then $(\bar{\bw}^{(k)},\bar{\bz}^{(k)})\in V(\tilde{I})$. Moreover, the two solutions possess the same multiplicity.
\end{lemma}
\begin{proof}
Using Lemma~\ref{associated} with $\lambda=0$, if $p_k(\bz^{(k)},\bw^{(k)})=\bar{p}_k(\bw^{(k)},\bz^{(k)})=0$, then $p_k(\bar{\bw}^{(k)},\bar{\bz}^{(k)})=\bar{p}_k(\bar{\bz}^{(k)},\bar{\bw}^{(k)})=0$.
\end{proof}

\begin{Def}
We call a pair of distinct points $(\bz,\bw),(\bar{\bw},\bar{\bz})\in V(\tilde{I})$ an \emph{associated pair}. We call a point $(\bz,\bar{\bz})\in V(\tilde{I})$ a \emph{conjugated single}.
\end{Def}

The following proposition links the problem to its real counterpart.

\begin{prop}\label{swap}
The endomorphism of $\C^{2r}$ given by
\begin{equation*}
(\bz,\bw)\mapsto\biggl(\frac{\bz+\bw}{2},\frac{\bz-\bw}{2i}\biggr)\qquad\text{with inverse}\qquad (\bx,\by)\mapsto(\bx+i\by,\bx-i\by),
\end{equation*}
sends the points in $V(\tilde{I})$ to the points in 
\begin{equation*}
    V(p_1^{\Ree}(\bx,\by),p_1^{\Imm}(\bx,\by),\ldots,p_c^{\Ree}(\bx,\by),p_c^{\Imm}(\bx,\by))
\end{equation*}
considering coordinates $\bz=\bx+i\by, \bw=\bx-i\by$, sending conjugated singles to real solutions and associated pairs to pairs of conjugated solutions.
\end{prop}
\begin{proof}
The thesis follows unfolding the definitions of the zero loci.
\end{proof}

We define a linear extension of the conjugation map
\begin{equation*}
    \ast\colon\C[\bz,\bw]\to\C[\bz,\bw],\qquad\ast(f)\mapsto f^{\ast}
\end{equation*}
such that $z_k^{\ast}=w_k$, $w_k^{\ast}=z_k$ for $k=1,\ldots,r$ and $a^{\ast}=\bar{a}$ for $a\in\C$. Clearly, it holds $\ast^2=\Id$ and from our choice of the ideal $\tilde{I}$, this map passes to the quotient $\C[\bz,\bw]/\tilde{I}$.

We will often consider polynomials $\xi=\xi^{\ast}\in\C[\bz,\bw]$. We prove a useful lemma.

\begin{lemma}
It holds $\xi=\xi^{\ast}\in\C[\bz,\bw]$ iff $\xi(\bx+i\by,\bx-i\by)\in\R[\bx,\by]$.
\end{lemma}
\begin{proof}
The thesis follows from the equalities
\begin{align*}
    \xi(\bx+i\by,\bx-i\by)&=\xi(\bz,\bw)=\xi(\bz,\bw)^{\ast}=\bar{\xi}(\bw,\bz)\\
    &=\bar{\xi}(\bx-i\by,\bx+i\by)=\overline{\xi(\bx+i\by,\bx-i\by)}.
\end{align*}
\end{proof}

We will simply use the notation $\xi(\bx,\by)$ to mean $\xi(\bx+i\by,\bx-i\by)$.

\begin{Def}
We define for $\xi=\xi^{\ast}\in\C[\bz,\bw]$ the sesquilinear forms
\begin{equation*}
  \C[\bz,\bw]/\tilde{I}\times\C[\bz,\bw]/\tilde{I}\to\C,\qquad ([f],[g])\mapsto\mathcal{K}_{\C}^{\xi}(f,g)\coloneqq\Tr(M_{\xi}M_fM_{g^{\ast}}).
\end{equation*}
We call the form $\mathcal{K}_{\C}^{1}$ the \emph{Hermitian Killing form}.  
\end{Def}

\begin{prop}\label{conjtr}
The form $\mathcal{K}_{\C}^{h}$ is Hermitian.
\end{prop}
\begin{proof}
The form is clearly sesquilinear.

We only need to show that for any $[f],[g]\in\C[\bz,\bw]/\tilde{I}$, it holds
\begin{equation*}
    \mathcal{K}_{\C}^{\xi}(f,g)=\Tr(M_{\xi fg^{\ast}})=\overline{\Tr(M_{\xi gf^{\ast}})}=\overline{\mathcal{K}_{\C}^{\xi}(g,f)}.
\end{equation*}
The second equality follows from
\begin{align*}
    \Tr(M_{\xi fg^{\ast}})&=\sum_{k=1}^{d}m_k\xi(\bz^{(k)},\bw^{(k)})f(\bz^{(k)},\bw^{(k)})\bar{g}(\bw^{(k)},\bz^{(k)})\tag{Lemma~\ref{prodeigen}, v}\\
    &=\sum_{k=1}^{d}m_k\overline{\bar{\xi}(\bar{\bz}^{(k)},\bar{\bw}^{(k)})\bar{f}(\bar{\bz}^{(k)},\bar{\bw}^{(k)})g(\bar{\bw}^{(k)},\bar{\bz}^{(k)})}\tag{Conjugation}\\
    &=\sum_{k=1}^{d}m_k\overline{\bar{\xi}(\bw^{(k)},\bz^{(k)})\bar{f}(\bw^{(k)},\bz^{(k)})g(\bz^{(k)},\bw^{(k)})}\tag{Lemma~\ref{asso}}\\
    &=\sum_{k=1}^{d}m_k\overline{\xi(\bz^{(k)},\bw^{(k)})\bar{f}(\bw^{(k)},\bz^{(k)})g(\bz^{(k)},\bw^{(k)})}\tag{$\xi=\xi^{\ast}$}\\
    &=\overline{\Tr(M_{\xi gf^{\ast}})}\tag{Lemma~\ref{prodeigen}, v}.
\end{align*}
\end{proof}

Consider the decomposition of Proposition~\ref{decA} with $\xi=\xi^{\ast}\in\C[\bz,\bw]$. 

\begin{lemma}\label{decA2}
The map $\ast$ fixes the sub algebras $A_k$ corresponding to a conjugated single and switches two sub algebras $A_k$ and $A_j$ corresponding to an associated pair. In other terms, it holds $e_k^{\ast}=e_k$ in the first case and $e_k^{\ast}=e_j$, $e_j^{\ast}=e_k$ in the latter.

Moreover, it holds $e_ke_k^{\ast}=e_k$ in the first case and $e_ke_k^{\ast}=0$ in the latter.
\end{lemma}
\begin{proof}
From the equalities 
\begin{equation*}
    \xi-\xi(\bz^{(k)},\bar{\bz}^{(k)})=\xi^{\ast}-\bar{\xi}(\bar{\bz}^{(k)},\bz^{(k)})=(\xi-\xi(\bz^{(k)},\bar{\bz}^{(k)}))^{\ast}
\end{equation*}
the part about conjugated singles follows.

From the equalities 
\begin{equation*}
    \xi-\xi(\bz^{(k)},\bw^{(k)})=\xi^{\ast}-\bar{\xi}(\bw^{(k)},\bz^{(k)})=(\xi-\xi(\bar{\bw}^{(k)},\bar{\bz}^{(k)}))^{\ast}
\end{equation*}
the part about associated pairs follows.

The last statement follows from Proposition~\ref{decA}.
\end{proof}

Let $d_1\leq d$ be the number of conjugated singles. We reorder the decomposition in such a way that
\begin{equation}\label{dec}
\C[\bz,\bw]/\tilde{I}=(\oplus_{k=1}^{d_1}A_k)\oplus(\oplus_{k=d_1+1}^{\frac{d-{d_1}}{2}}(A_k\oplus A_k^{\ast}))
\end{equation}
where the first $d_1$ summands are the sub algebras corresponding to conjugated singles and the second summands correspond to associated pairs.

\begin{prop}\label{decompo}
For $\xi=\xi^{\ast}\in\C[\bz,\bw]$, the decomposition \eqref{dec} is orthogonal with respect to $\mathcal{K}_{\C}^{\xi}$.
\end{prop}
\begin{proof}
The thesis follows since the element of the decomposition are invariant under the action of the map $\ast$, and $\mathcal{K}_{\C}^{\xi}$ behaves as $\mathcal{K}_{\R}^{\xi}$ on this decomposition.

However, we prove this statement directly. Computing the form on the units is sufficient, we consider three cases:
\begin{enumerate}[\normalfont i)]
    \item if $k\neq j\leq d_1$, then $\mathcal{K}_{\C}^{\xi}(e_k,e_j)=\Tr(M_{\xi e_ke_j^{\ast}})=\Tr(M_{\xi e_ke_j})=\Tr(M_{0})=0$,
    \item if $k\leq d_1$ and $j>d_1$, then $\mathcal{K}_{\C}^{\xi}(e_k,e_j+e_j^{\ast})=\Tr(M_{\xi(e_ke_j^{\ast}+e_ke_j)})=\Tr(M_{0})=0$,
    \item if $k\neq j>d_1$, then $\mathcal{K}_{\C}^{\xi}(e_k+e_k^{\ast},e_j+e_j^{\ast})=\Tr(M_{\xi(e_k+e_k^{\ast})(e_j^{\ast}+e_j)})=\Tr(M_{0})=0$.
\end{enumerate}
These equalities follow from the properties of the units $e_k$ for $k=1,\ldots,d$ exploited in Lemma~\ref{decA2}.
\end{proof}

The following is the main result of this work, its proof is analogous to the proof of Theorem~\ref{hermite}.

\begin{teo}\label{comphermite}
Let $\xi=\xi^{\ast}\in\C[\bz,\bw]$, then:
\begin{enumerate}[\normalfont i)]
    \item The variety  $V(\tilde{I})$ consists only in conjugated singles $(\bz^{(k)},\bar{\bz}^{(k)})$ with multiplicity $1$ such that $\xi(\bz^{(k)},\bar{\bz}^{(k)})>0$ iff $\mathcal{K}_{\C}^{\xi}$ is positive definite.
    \item The rank of $\mathcal{K}_{\C}^{\xi}$ is the number of distinct points $(\bz^{(k)},\bw^{(k)})\in V(\tilde{I})$ such that $\xi(\bz^{(k)},\bw^{(k)})\neq 0$.
    \item The number of distinct conjugated singles $(\bz^{(k)},\bar{\bz}^{(k)})\in V(\tilde{I})$ such that $\xi(\bz^{(k)},\bar{\bz}^{(k)})>0$ minus the number of distinct conjugated singles $(\bz^{(k)},\bar{\bz}^{(k)})\in V(\tilde{I})$ such that $\xi(\bz^{(k)},\bar{\bz}^{(k)})<0$ equals the signature of $\mathcal{K}_{\C}^{\xi}$.
\end{enumerate}
If moreover $\xi$ is such that the evaluations $\xi(\bz^{(k)},\bw^{(k)})$ are all different for $k=1,\ldots,d$, then:
\begin{enumerate}[\normalfont i)]
    \setcounter{enumi}{3}
    \item The number of distinct conjugated singles in $(\bz^{(k)},\bar{\bz}^{(k)})\in V(\tilde{I})$ such that $\xi(\bz^{(k)},\bar{\bz}^{(k)})>0$ equals the number of positive eigenvalues of $\mathcal{K}_{\C}^{\xi}$ minus the number of negative eigenvalues of $\mathcal{K}_{\C}^1$.
    \item The number of distinct conjugated singles in $(\bz^{(k)},\bar{\bz}^{(k)})\in V(\tilde{I})$ such that $\xi(\bz^{(k)},\bar{\bz}^{(k)})<0$ equals the number of negative eigenvalues of $\mathcal{K}_{\C}^{\xi}$ minus the number of negative eigenvalues of $\mathcal{K}_{\C}^1$.
\end{enumerate}
\end{teo}
\begin{proof}
From Proposition~\ref{decompo} we can study the signature of $\mathcal{K}_{\C}^{\xi}$ over the decomposition \eqref{dec}.

For $k=1,\ldots,d$, consider a basis of $A_k$ given by $\lbrace [e_kf_{k,j}]\rbrace_{j=0,\ldots,m_{k}-1}$, where $f_{k,j}\in\C[\bz,\bw]$ are polynomials such that $f_{k,0}\equiv 1$ and $f_{k,j}(\bz^{(k)},\bw^{(k)})=0$ for $j=1,\ldots,m_k-1$.

If $k\leq d_1$, then $f_{k,j}(\bz^{(k)},\bar{\bz}^{(k)})^{\ast}=0$ and thus the only non vanishing term is 
\begin{align*}
    \Tr(M_{\xi e_ke_k^{\ast}})=\Tr(M_{\xi e_k^2})=\Tr(M_{\xi e_k})&=\sum_{j=1}^dm_j\xi(\bz^{(j)},\bw^{(j)})\delta_{k,j}\\
    &=m_k\xi(\bz^{(k)},\bar{\bz}^{(k)})
\end{align*}
corresponding to the leading term.

If $k>d_1$, consider the basis of the subspace $A_k\oplus A_{k}^{\ast}$ given by 
\begin{equation*}
\lbrace [e_{k}f_{k,j}+(e_{k}f_{k,j})^{\ast}],[(e_kf_{k,j}-(e_kf_{k,j})^{\ast})/i]\rbrace_{j=0,\ldots,m_{k}-1}.
\end{equation*}
Similarly as before, we can consider just the subspace spanned by $\lbrace[e_k+e_k^{\ast}],[(e_k-e_k^{\ast})/i]\rbrace$ corresponding to the second leading principal minor that is
\begin{align*}
  &\begin{bmatrix}
      \Tr(M_{\xi(e_k+e_k^{\ast})(e_k^{\ast}+e_k)}) & i\Tr(M_{\xi(e_k+e_k^{\ast})(e_k^{\ast}-e_k)})\\
      -i\Tr(M_{\xi(e_k-e_k^{\ast})(e_k^{\ast}+e_k)}) & \Tr(M_{\xi(e_k-e_k^{\ast})(e_k^{\ast}-e_k)})
  \end{bmatrix}=\begin{bmatrix}
      \Tr(M_{\xi(e_k+e_k^{\ast})}) & -i\Tr(M_{\xi(e_k-e_k^{\ast})})\\
      -i\Tr(M_{\xi(e_k-e_k^{\ast})}) & -\Tr(M_{\xi(e_k+e_k^{\ast})})
  \end{bmatrix}
\end{align*}
and from $e_k(\bz^{(j)},\bw^{(j)})=\delta_{k,j}$ and  Lemma~\ref{asso}, this matrix is
\begin{equation*}
    m_k\begin{bmatrix}
      \xi(\bz^{(k)},\bw^{(k)})+\xi(\bar{\bw}^{(k)},\bar{\bz}^{(k)}) & -i(\xi(\bz^{(k)},\bw^{(k)})-\xi(\bar{\bw}^{(k)},\bar{\bz}^{(k)}))\\
      -i(\xi(\bz^{(k)},\bw^{(k)})-\xi(\bar{\bw}^{(k)},\bar{\bz}^{(k)})) & -(\xi(\bz^{(k)},\bw^{(k)})+\xi(\bar{\bw}^{(k)},\bar{\bz}^{(k)}))
  \end{bmatrix}.
\end{equation*}
Since from $\xi=\xi^{\ast}$ the equality $\xi(\bz^{(k)},\bw^{(k)})=\overline{\xi(\bar{\bw}^{(k)},\bar{\bz}^{(k)})}$ holds, we rewrite this matrix 
as
\begin{equation*}
    m_k\begin{bmatrix}
      \xi(\bz^{(k)},\bw^{(k)})+\overline{\xi(\bz^{(k)},\bw^{(k)})} & -i(\xi(\bz^{(k)},\bw^{(k)})-\overline{\xi(\bz^{(k)},\bw^{(k)})})\\
      -i(\xi(\bz^{(k)},\bw^{(k)})-\overline{\xi(\bz^{(k)},\bw^{(k)})}) & -(\xi(\bz^{(k)},\bw^{(k)})+\overline{\xi(\bz^{(k)},\bw^{(k)})})
  \end{bmatrix}
\end{equation*}
and note that it is real symmetric. We decompose it as
\begin{equation*}
  m_kU\begin{bmatrix}
      \xi(\bz^{(k)},\bw^{(k)}) & 0\\
      0 & \overline{\xi(\bz^{(k)},\bw^{(k)})}
  \end{bmatrix}U^{T}
\end{equation*}
where $U=\begin{bmatrix}
      1 & 1\\
      -i & i
  \end{bmatrix}$. Using the Binet theorem the determinant is
\begin{equation*}
    (2i)^2m_k^2|\xi(\bz^{(k)},\bw^{(k)})|^2=-4m_k^2|\xi(\bz^{(k)},\bw^{(k)})|^2
\end{equation*}
and the thesis follows.
\end{proof}

\begin{coro}\label{count}
The number of zeros of a generalized polynomial system $p_k(\bz,\bar{\bz})=0$ for $k=1,\ldots,c$, is equal to the number of positive minus the number of negative eigenvalues of the Hermitian Killing form $\mathcal{K}_{\C}^1$ of the system
\begin{equation*}
\begin{cases}
    p_k(\bz,\bw)=0 &k=1,\ldots c\\
    p_k(\bz,\bw)^{\ast}=\bar{p}_k(\bw,\bz)=0 &k=1,\ldots c
\end{cases}.
\end{equation*}
Moreover, the signature of $\mathcal{K}_{\C}^1$ is equal to the signature of the Killing form $\mathcal{K}_{\R}^1$ considering coordinates $\bz=\bx+i\by$, $\bw=\bx-i\by$.
\end{coro}
\begin{proof}
The thesis follows from Theorem~\ref{comphermite} and Proposition~\ref{swap}.
\end{proof}

\vspace{1cm}
We want to compare the Hermitian Killing form and the Killing form.

\begin{ex}\label{ex1}
Consider the case of a second degree harmonic equation
\begin{equation}\label{prob}
h(z,\bar{z})=z^2+a\bar{z}+b=0\qquad\text{for $a,b\in\C$.}
\end{equation}
Denoting the real and imaginary parts of the coefficients with the subscripts $1$ and $2$ respectively, equation \eqref{prob} leads us to the real system
\begin{equation}\label{realsys}
\begin{cases}
    h^{\Ree}(x,y)=x^2-y^2+a_1x+a_2y+b_1=0\\
    h^{\Imm}(x,y)=2xy+a_2x-a_1y+b_2=0
\end{cases}.
\end{equation}
We compute the matrix representing the Killing form associated to the real system \eqref{realsys} with respect to the basis $\lbrace [1],[x],[y],[y^2]\rbrace$ of $\R[x,y]/\langle h^{\Ree},h^{\Imm}\rangle$ to obtain the real symmetric matrix
\begin{equation*}
\mathcal{K}_{\R}^1=\begin{bmatrix}
    4 & 0 & 0 & \frac{3|a|^2+4b_1}{2}\\
    0 & \frac{3|a|^2-4b_1}{2} & -2b_2 & \frac{3|a|^2a_1+4(a_1b_1+a_2b_2)}{4}\\
    0 & -2b_2 & \frac{3|a|^2+4b_1}{2} & \frac{3|a|^2a_2+12(a_2b_1-a_1b_2)}{4}\\
    \frac{3|a|^2+4b_1}{2} & \frac{3|a|^2a_1+4(a_1b_1+a_2b_2)}{4} & \frac{3|a|^2a_2+12(a_2b_1-a_1b_2)}{4} & \frac{9|a|^4+8(3a_1^2b_1+4a_2^2b_1-a_1a_2b_2+2b_1^2+bb_2^2)}{8}
\end{bmatrix}.
\end{equation*}
This matrix is quite articulate, even if it comes from an apparently simple problem.

On the other hand, the Hermitian matrix representing the Hermitian Killing form of the system
\begin{equation*}
\begin{cases}
    h(z,w)=z^2+aw+b=0\\
    h(z,w)^{\ast}=w^2+\bar{a}z+\bar{b}=0
\end{cases}
\end{equation*}
obtained by equation \eqref{prob}, with respect to the basis $\lbrace [1],[z],[w],[zw]\rbrace$ of $\C[z,w]/\langle h,h^{\ast}\rangle$ is
\begin{equation*}
\mathcal{K}_{\C}^{1}=\begin{bmatrix}
    4 & 0 & 0 & 3|a|^2\\
    0 & 3|a|^2 & -4b & 4a\bar{b}\\
    0 & -4\bar{b} & 3|a|^2 & 4\bar{a}b\\
    
    3|a|^2 & 4\bar{a}b & 4a\bar{b} & 3|a|^4+4|b|^2
\end{bmatrix}.
\end{equation*}
The signs of the second, third and fourth principal minors are given by 
\begin{enumerate}[\normalfont i)]
    \item $12|a|^2$,
    \item $3|a|^2-4|b|$,
    \item $27|a|^8-32|a|^4|b|^2-256|a^2\bar{b}+b^2|^2$,
\end{enumerate}
respectively. If any of them is equal to zero, equation \eqref{prob} admits at most $3$ solutions since the matrix posses at most $3$ positive eigenvalues. If any of them is negative, equation \eqref{prob} admits at most $3-1=2$ solutions since the matrix possesses at most $3$ positive eigenvalues and at least a negative eigenvalue. On the other hand, if i,ii and iii are all positive, which is the case for $a\neq 0$ and $b=0$, the Sylvester's Theorem assures that this matrix is positive definite and equation \eqref{prob} admits exactly $4$ solutions.   
\end{ex}

We can pass from the matrix representing $K_{\R}^{\xi}$ to the matrix representing $K_{\C}^{\xi}$ by applying the maps of Proposition~\ref{swap}, starting from a real symmetric matrix, this modification yields an Hermitian matrix and vice versa. Sometimes this procedure translates in a change of basis as exploited in the following proposition.

For the sake of simplicity, we will use a polynomial $f$ also denoting the value $\Tr(M_f)$, when used as entry in a matrix. Since it is irrelevant for what follows in this section, from now on we will drop the dependence from $\xi$ in the computations.

\begin{prop}\label{kckr}
For $n\in\N^{r}$, let $\C[\bz,\bw]/\tilde{I}=\lbrace[\bz^{\alpha}\bw^{\beta}]\mid\alpha_k+\beta_k\leq n_k\rbrace$ and $\xi=\xi^{\ast}\in\C[\bz,\bw]$, so that $\xi(\bx,\by)\in\R[\bx,\by]$. Then $\mathcal{K}^{\xi}_{\C}$ and $\mathcal{K}^{\xi}_{\R}$ can be obtained one from the other by a change of basis.

Moreover, this holds true for any restriction on $\lbrace[\bz^{\alpha}\bw^{\beta}]\mid\alpha_k+\beta_k\leq\tilde{n}_k\rbrace$ with $\tilde{n}_k\leq n_k$ for any $k=1,\ldots,r$.
\end{prop}
\begin{proof}
Using the map of Proposition~\ref{swap} and Lemma~\ref{prodeigen} we obtain 
\begin{align*}
    \bz^{\alpha}\bw^{\beta}=(\bx+i\by)^{\alpha}(\bx-i\by)^{\beta}=\sum_{\lbrace(\tilde{\alpha},\tilde{\beta})\mid\tilde{\alpha}+\tilde{\beta}=\alpha+\beta\rbrace}\gamma_{\tilde{\alpha},\tilde{\beta}}\bx^{\tilde{\alpha}}\by^{\tilde{\beta}}
\end{align*}
for suitables $\gamma_{\tilde{\alpha},\tilde{\beta}}$. Thus, if 
\begin{equation*}
    I=\langle p_1^{\Ree}(\bx,\by),p_1^{\Imm}(\bx,\by),\ldots,p_c^{\Ree}(\bx,\by),p_c^{\Imm}(\bx,\by)\rangle
\end{equation*}
then $\C[\bx,\by]/I=\lbrace[\bx^{\alpha}\by^{\beta}]\mid\alpha_k+\beta_k\leq n_k\rbrace$ and vice versa. We can obtain one basis from the other by linear combinations, if $\Gamma=(\gamma_{\tilde{\alpha},\tilde{\beta}})_{\lbrace(\tilde{\alpha},\tilde{\beta})\mid\tilde{\alpha}+\tilde{\beta}\leq n_k\rbrace}$ then
\begin{equation*}
    \Gamma\begin{bmatrix}
        \vdots\\
        \bx^{\alpha}\by^{\beta}\\
        \vdots
    \end{bmatrix}\begin{bmatrix}
        \cdots & \bx^{\alpha}\by^{\beta} & \cdots
    \end{bmatrix}\Gamma^H=\begin{bmatrix}
        \vdots\\
        \bz^{\alpha}\bw^{\beta}\\
        \vdots
    \end{bmatrix}\begin{bmatrix}
        \cdots & \bw^{\alpha}\bz^{\beta} & \cdots
    \end{bmatrix}.
\end{equation*}
\end{proof}

\begin{ex}\label{change}
We consider the simplest non trivial case of Proposition~\ref{kckr}. Thus, let $r=1$ and $n=1$, they lead us to the bases $\lbrace [1],[z],[w]\rbrace$ and $\lbrace [1],[x],[y]\rbrace$. Using the notation of $U$ from the proof of Theorem~\ref{comphermite} we obtain
\begin{align*}
\begin{bmatrix}
    1 & \begin{matrix} 0 & 0 \end{matrix}\\
    \begin{matrix} 0 \\ 0 \end{matrix} & U^H\\
\end{bmatrix}\begin{bmatrix}
    1 & x & y\\
    x & x^2 & xy\\
    y & xy & y^2
\end{bmatrix}\begin{bmatrix}
    1 & \begin{matrix} 0 & 0 \end{matrix}\\
    \begin{matrix} 0 \\ 0 \end{matrix} & U\\
\end{bmatrix}&=\begin{bmatrix}
    1 & 0 & 0\\
    0 & 1 & i\\
    0 & 1 & -i
\end{bmatrix}\begin{bmatrix}
    1 & x & y\\
    x & x^2 & xy\\
    y & xy & y^2
\end{bmatrix}\begin{bmatrix}
    1 & 0 & 0\\
    0 & 1 & 1\\
    0 & -i & i
\end{bmatrix}\\
&=\begin{bmatrix}
    1 & x-iy & x+iy\\
    x+iy & x^2+y^2 & x^2-y^2+i2xy\\
    x-iy & x^2-y^2-i2xy & x^2+y^2
\end{bmatrix}
\end{align*}
and
\begin{align*}
\begin{bmatrix}
    1 & \begin{matrix} 0 & 0 \end{matrix}\\
    \begin{matrix} 0 \\ 0 \end{matrix} & \frac{1}{2}U\\
\end{bmatrix}\begin{bmatrix}
    1 & w & z\\
    z & zw & z^2\\
    w & w^2 & zw
\end{bmatrix}\begin{bmatrix}
    1 & \begin{matrix} 0 & 0 \end{matrix}\\
    \begin{matrix} 0 \\ 0 \end{matrix} & \frac{1}{2}U^H\\
\end{bmatrix}&=\begin{bmatrix}
    1 & 0 & 0\\
    0 & \frac{1}{2} & \frac{1}{2}\\
    0 & \frac{1}{2i} & -\frac{1}{2i}
\end{bmatrix}\begin{bmatrix}
    1 & w & z\\
    z & zw & z^2\\
    w & w^2 & zw
\end{bmatrix}\begin{bmatrix}
    1 & 0 & 0\\
    0 & \frac{1}{2} & -\frac{1}{2i}\\
    0 & \frac{1}{2} & \frac{1}{2i}
\end{bmatrix}\\
&=\begin{bmatrix}
    1 & \frac{w+z}{2} & \frac{z-w}{2i}\\
    \frac{w+z}{2} & \frac{2zw+z^2+w^2}{4} & \frac{z^2-w^2}{4i}\\
    \frac{z-w}{2i} & \frac{z^2-w^2}{4i} & \frac{2zw-z^2-w^2}{4}
\end{bmatrix}.
\end{align*}   
\end{ex}
The change of basis of this last example is a special case, which is particularly simple and will be valid whenever we restrict the form on these subspaces. However, in general such a change of basis is not trivial, as in the case of Example~\ref{ex1} with bases $\lbrace [1],[z],[w],[zw]\rbrace$ and $\lbrace [1],[x],[y],[y^2]\rbrace$.\vspace{1cm}

Note that since for any polynomial $f$ holds $\Tr(M_f)=\overline{\Tr(M_{f^{\ast}})}$, when computing the matrix representing $\mathcal{K}_{\C}^{\xi}$ the entries above the diagonal could be conjugated each other. In Example~\ref{ex1} we have checked that $\mathcal{K}_{\C}^1$ is simpler to be computed than $\mathcal{K}_{\R}^1$. This check extends to all the cases we have computed and it seems that the form $K_{\C}^{\xi}$, beyond its theoretical use in the proof of Theorem~\ref{main}, has also a computational advantage.

\begin{ex}\label{ex2}
For the bases $\lbrace [1],[z],[w],[zw]\rbrace$ and $\lbrace [1],[z],[w],[z^2]\rbrace$ with matrices
\begin{equation*}
\begin{bmatrix}
    1 & w & \underline{z} & \underline{zw}\\
    z & zw & \underline{z^2} & \underline{z^2w}\\
    w & w^2 & zw & zw^2\\
    zw & zw^2 & z^2w & \underline{z^2w^2}
\end{bmatrix}\qquad\text{and}\qquad\begin{bmatrix}
    1 & w & \underline{z} & \underline{z^2}\\
    z & \underline{zw} & z^2 & \underline{z^3}\\
    w & w^2 & zw & \underline{z^2w}\\
    w^2 & w^3 & zw^2 & \underline{z^2w^2}
\end{bmatrix},
\end{equation*}
we underline a set of $5$ and $6$ entries respectively, which are generically different and non conjugated.

For the bases $\lbrace [1],[x],[y],[y^2]\rbrace$ and $\lbrace [1],[x],[y],[xy]\rbrace$ with matrices
\begin{equation*}
\begin{bmatrix}
    1 & \underline{x} & \underline{y} & \underline{y^2}\\
    x & \underline{x^2} & \underline{xy} & \underline{xy^2}\\
    y & xy & y^2 & \underline{y^3}\\
    y^2 & xy^2 & y^3 & \underline{y^4}
\end{bmatrix}\qquad\text{and}\qquad\begin{bmatrix}
    1 & \underline{x} & \underline{y} & \underline{xy}\\
    x & \underline{x^2} & xy & \underline{x^2y}\\
    y & xy & \underline{y^2} & \underline{xy^2}\\
    xy & x^2y & xy^2 & \underline{x^2y^2}
\end{bmatrix},
\end{equation*}
we underline a set of $8$ entries each, which are generically different and non conjugated.

If we consider the two left matrices, we are in the case of Example~\ref{ex1}.

Note that we have more control in the diagonal. In fact, the entries in the diagonal of the matrix representing $K_{\R}^{\xi}$ are different in general, for two monomials $\bx^{\alpha_1}\by^{\beta_1}\neq \bx^{\alpha_2}\by^{\beta_2}$ it holds $(\bx^{\alpha_1}\by^{\beta_2})^2\neq(\bx^{\alpha_2}\by^{\beta_2})^2$. On the other hand, in the matrix representing $K_{\C}^{\xi}$, it happens to have equal entries in the diagonal, for two monomials $\bz^{\alpha_1}\bw^{\beta_1}\neq \bz^{\alpha_2}\bw^{\beta_2}$ the equality $\bz^{\alpha_1}\bw^{\beta_1}\cdot( \bz^{\alpha_1}\bw^{\beta_1})^{\ast}=\bz^{\alpha_2}\bw^{\beta_2}\cdot( \bz^{\alpha_2}\bw^{\beta_2})^{\ast}$ holds true when $\alpha_1+\beta_1=\alpha_2+\beta_2$.
\end{ex}

%-------------------------------------------------------------------------------------------

\section{Harmonic polynomial equations}\label{harm}

We set $r=1$ and consider harmonic polynomial equations $h(z,\bar{z})=q(z)+p(\bar{z})=0$, where
\begin{equation}\label{harmpol2}
 q(z)=b_mz^m+\ldots+b_1z\qquad\text{and}\qquad p(\bar{z})=\bar{z}^n+a_{n-1}\bar{z}^{n-1}+\ldots+a_1\bar{z}+a_0 
\end{equation}
are complex polynomials. The ideal $\tilde{I}=\langle h(z,w),h(z,w)^{\ast}\rangle$ is generically zero dimensional. Suppose $m=\deg(q)\leq\deg(p)=n$, then the collection $\lbrace [z^{\alpha}w^{\beta}]\mid\alpha,\beta< n\rbrace$ of cardinality $n^2$ is a basis of the quotient space $\C[z,w]/\tilde{I}$. We want to bound the number of conjugated singles. 

We sharpen the trivial bound $n^2$ with the following result.

\begin{teo}\label{bound}
Using the notations of equations \eqref{harmpol2}, if $n-2\geq m$ the harmonic equation $h(z,\bar{z})=q(z)+p(\bar{z})=0$ admits at most $n^2-1$ solutions if $(n-1)a_{n-1}^2-2na_{n-2}=0$ and at most $n^2-2$ solutions otherwise.
\end{teo}

We deal with two lemmas before proving Theorem~\ref{bound}.

\begin{lemma}\label{lem1}
If $n\geq 2$ with $n-1\geq m$, for the harmonic equation $h(z,\bar{z})=q(z)+p(\bar{z})=0$ hold 
\begin{equation*}
\Tr(M_w)=-na_{n-1}\qquad\text{and}\qquad\Tr(M_{w^2})=n(a_{n-1}^2-2a_{n-2}).
\end{equation*}
\end{lemma}
\begin{proof}
Consider the basis $\lbrace [z^{\alpha}w^{\beta}]\mid\alpha,\beta< n\rbrace$.

Compute $\Tr(M_w)$:

The multiplications $w\cdot z^{k}w^{j}$ with $0\leq k\leq n-1$ and $0\leq j\leq n-2$ do not give contribution.

Let $0\leq k\leq n-1$, in $\C[z,w]/\tilde{I}$ hold the equalities
\begin{align*}
  w\cdot z^{k}w^{n-1}&=z^kw^{n}=z^k
  (w^n-h(z,w))=-\sum_{j=1}^m{b_{j}}z^{j+k}-\sum_{s=0}^{n-1}a_sz^kw^{s}\\
  &=-\sum_{j=k+1}^{n-1}b_{j-k}z^{j}-\sum_{\ell=n}^{m+k}b_{\ell-k}z^{\ell}-\sum_{s=0}^{n-1}a_sz^{k}w^{s}.
\end{align*}
The second sum appears if $n\leq m+k$, in this case the element $z^{\ell}$ with $n\leq \ell\leq m+k\leq 2n-2$ can be rewritten as $z^{\ell-n}(z^n-\bar{h}(w,z))$, and the higher degree of $z$ for monomials with variable $w$ is $\ell-n\leq m+k-n<k$. In the end, for any $k$, the term $-a_{n-1}z^kw^{n-1}$ is given by the third sum and then 
\begin{equation*}
    \Tr(M_w)=\sum_{k=0}^{n-1}(-a_{n-1})=-na_{n-1}.
\end{equation*}

Compute $\Tr(M_{w^2})$:

The multiplications $w^2\cdot z^{k}w^{j}$ with $0\leq k\leq n-1$ and $0\leq j\leq n-3$ do not give contribution.
Let $0\leq k\leq n-1$, for $w^2\cdot z^kw^{n-2}=z^kw^n$ we consider the computations for $\Tr(M_{w})$ and get the term $-a_{n-2}z^kw^{n-2}$. Now, in $\C[z,w]/\tilde{I}$ hold the equalities

\begin{align*}
  w^2\cdot z^{k}w^{n-1}&=z^kw^{n+1}=z^kw(w^n-h(z,w))=-\sum_{j=1}^mb_jz^{j+k}w-\sum_{s=0}^{n-1}a_sz^kw^{s+1}\\
  &=-\sum_{j=k+1}^{n-1}b_{j-k}z^{j}w-\sum_{\ell=n}^{m+k}b_{\ell-k}z^{\ell}w-\sum_{s=0}^{n-2}a_sz^{k}w^{s+1}-a_{n-1}z^k(w^n-h(z,w)).
\end{align*}
The second sum appears if $n\leq m+k$, again we rewrite $z^{\ell}$ and the higher degree of $z$ for monomials with variable $w$ is $\ell-n\leq m+k-n<k$. Thus, adding the coefficients obtained before to the coefficients given by the third and fourth sum we get 
\begin{equation*}
\Tr(M_{w^2})=\sum_{k=0}^{n-1}(-a_{n-2})+\sum_{j=0}^{n-1}\left(-a_{n-2}+a_{n-1}^2\right)=n(a_{n-1}^2-2a_{n-2}).
\end{equation*}
\end{proof}

\begin{lemma}\label{lem2}
If $n\geq 2$, for the harmonic equation $h(z,\bar{z})=q(z)+p(\bar{z})=0$ hold
\begin{equation*}
\Tr(M_{zw})=\begin{cases}
    |a_{n-1}|^2+(2n-1)|b_{n-1}|^2&\text{if $n-1=m$}\\
    |a_{n-1}|^2&\text{if $n-2\geq m$}\\
\end{cases}.
\end{equation*}
\end{lemma}
\begin{proof}
Consider the basis $\lbrace [z^{\alpha}w^{\beta}]\mid\alpha,\beta< n\rbrace$.

Compute $\Tr(M_{zw})$:

The multiplications $zw\cdot z^{k}w^{j}$ with $0\leq k\leq n-2$ and $0\leq j\leq n-2$ do not give contribution.

Let $1\leq k\leq n-1$, for $zw\cdot z^{k-1}w^{n-1}=z^{k}w^{n}$, we consider the computations for $\Tr(M_{w})$ in Lemma~\ref{lem1}. Thus, if $m=n-1$, we get the term $|b_{n-1}|^2z^{k-1}w^{n-1}$ given by the monomial $b_{m}z^{m+k}=b_{m}z^{n+k-1}$ in the second sum, otherwise we get no terms. By conjugation we get the information on the products $zw\cdot z^{n-1}w^{k-1}=z^{n}w^{k}$ with $1\leq k\leq n-1$. Lastly, in $\C[z,w]/\tilde{I}$ hold the equalities
\begin{align*}
  zw\cdot z^{n-1}w^{n-1}&=z^nw^n=(z^n-\bar{h}(w,z))(w^n-h(z,w))\\
  &=\left(\sum_{j=1}^m\bar{b}_jw^{j}+\sum_{t=0}^{n-1}\bar{a}_{t}z^{t}\right)\left(\sum_{k=1}^mb_kz^{k}+\sum_{s=0}^{n-1}a_sw^{s}\right)\\
  &=\sum_{j,k=1}^mb_k\bar{b}_jz^{k}w^{j}+\sum_{\ell=1}^{n+m-1}(\gamma_{\ell}z^{\ell}+\bar{\gamma}_{\ell}w^{\ell})+\sum_{s,t=0}^{n-1}\bar{a}_ta_sz^{t}w^{s},
\end{align*}
for suitable $\gamma_{\ell}$. The first sum yields the term $|b_{n-1}|^2z^{n-1}w^{n-1}$ if and only if $m=n-1$. The last sum yields the term $|a_{n-1}|^2z^{n-1}w^{n-1}$. Arguing as in Lemma~\ref{lem1}, we rewrite the elements $z^{\ell}$ with $\ell\geq n$ and the higher degree of $z$ for monomials with variable $w$ is $\ell-n\leq n+m-1-n=m-1<n-1$, similarly for the terms $w^{\ell}$. In the end, adding the coefficients, if $m<n-1$ we get 
\begin{equation*}
    \Tr(M_{zw})=|a_{n-1}|^2,
\end{equation*}
while if $m=n-1$ we get
\begin{equation*}
  \Tr(M_{zw})=|b_{n-1}|^2+2\sum_{k=1}^{n-1}|b_{n-1}|^2+|a_{n-1}|^2=|a_{n-1}|^2+(2n-1)|b_{n-1}|^2.
\end{equation*}
\end{proof}

Now, the proof of Theorem~\ref{bound}.

\begin{proof}
We consider the matrix representing the Hermitian Killing form on the subspace spanned by $\lbrace[1],[z],[w]\rbrace$. From the choice of the basis follows that to compute this matrix it is sufficient to calculate $\Tr(M_w),\Tr(M_{w^2})$ and $\Tr(M_{zw})$. Then, by Lemma~\ref{lem1} and \ref{lem2} the matrix is
\begin{equation*}
\begin{bmatrix}
    \Tr(M_1) & \Tr(M_w) & \Tr(M_z)\\
    \Tr(M_z) & \Tr(M_{zw}) & \Tr(M_{z^2})\\
   \Tr(M_w) & \Tr(M_{w^2}) & \Tr(M_{zw})
\end{bmatrix}=\begin{bmatrix}
    n^2 & -na_{n-1} & -n\bar{a}_{n-1}\\
    -n\bar{a}_{n-1} & |a_{n-1}|^2 & n(\bar{a}_{n-1}^2-2\bar{a}_{n-2})\\
    -na_{n-1} & n(a_{n-1}^2-2a_{n-2}) & |a_{n-1}|^2
\end{bmatrix}
\end{equation*}
and its determinant is $-n^2|(n-1)a_{n-1}^2-2na_{n-2}|^2$. In particular, the matrix possesses a negative eigenvalue if the determinant does not vanish and a non positive eigenvalue otherwise. Thus, by Corollary~\ref{count}, the equation admits at most $(n^2-1)-1=n^2-2$ solutions in the first case and $n^2-1$ in the latter.
\end{proof}

\begin{rmk}
For $n=2$ the formula $(n-1)a_{n-1}^2-2na_{n-2}=0$ in Theorem~\ref{bound} is the discriminant $a_{n-1}^2-4a_{n-2}$ of the second degree polynomial $z^2+a_{n-1}z+a_{n-2}$. For general $n\geq 2$, the formula is a multiple of the discriminant of the second degree polynomial given by the $(n-2)$-th derivative. In fact
\begin{align*}
    \frac{\partial^{n-2}z^n+a_{n-1}z^{n-1}+a_{n-2}z^{n-2}+\ldots+a_0}{\partial z^{n-2}}&=\frac{n!}{2}z^2+(n-1)!a_{n-1}z+(n-2)!a_{n-2}\\
    &=(n-2)!\left(\frac{n(n-1)}{2}z^2+(n-1)a_{n-1}z+a_{n-2}\right)
\end{align*}
and the discriminant of the polynomial in the parenthesis is 
\begin{align*}
  (n-1)^2a_{n-1}^2-2n(n-1)a_{n-2}=(n-1)((n-1)a_{n-2}^2-2na_{n-2}).
\end{align*}

Another natural appearance of this formula is in the third term of the logarithmic discriminant of a $n$ dimensional vector bundle $\mathcal{E}$. From the equalities
\begin{equation*}
    \ch(\mathcal{E})=n+c_1+\frac{c_1^2-2c_2}{2}+\ldots=n\left(1+\frac{c_1}{n}+\frac{c_1^2-2c_2}{2n}+\ldots\right)
\end{equation*}
and the Maclaurin series $\log(1+z)=z-z^2/2+\ldots$ we get 
\begin{align*}
    \log(\ch(\mathcal{E}))&=\log n+\log\left(1+\frac{c_1}{n}+\frac{c_1^2-2c_2}{2n}+\ldots\right)\\
    &=\log n+\frac{c_1}{n}+\frac{c_1^2-2c_2}{2n}-\frac{c_1^2}{2n^2}+\ldots\\
    &=\log n+\frac{c_1}{n}+\frac{(n-1)c_1^2-2nc_2}{2n^2}+\ldots
\end{align*}

Note that the expression $(n-1)a_{n-1}^2-2na_{n-2}$ equals 
\begin{equation*}
    \frac{n^2\Tr(M_{w^2})-\Tr(M_{w})^2}{n^2}=\frac{n^2\Tr(M_{w}^2)-\Tr(M_{w})^2}{n^2}.
\end{equation*}
\end{rmk}

\subsection{A two parameter family with \textit{n=m}}

We now briefly consider the harmonic polynomial equation
\begin{equation*}
  h(z,\bar{z})=z^n+a\bar{z}^n+b=0\qquad\text{for $a,b\in\C$ with $|a|\neq 1$.}
\end{equation*}
With this choices the ideal $\langle h(z,w),h(z,w)^{\ast}\rangle\subset\C[z,w]$ is zero dimensional and the quotient $\C[z,w]/\langle h(z,w),h(z,w)^{\ast}\rangle$ possesses basis $\lbrace [z^{\alpha}w^{\beta}]\mid\alpha,\beta<n\rbrace$.

\begin{lemma}\label{deg1}
The harmonic equation $z+a\bar{z}+b=0$ for $a,b\in\C$ admits one solution iff $|a|\neq 1$, in this case if $b=0$ the solution is $z=0$.

If $|a|=1$ the equation admits an infinite amount of solutions if $b/\sqrt{a}\in\R$ and zero solutions otherwise.
\end{lemma}
\begin{proof}
Using real coordinates $z=x+iy$ and $a=a_1+ia_2, b=b_1+ib_2$, we divide the equation into its real and imaginary part to get the linear system
\begin{equation*}
\begin{bmatrix}
    1+a_1 & a_2\\
    a_2 & 1-a_1
\end{bmatrix}\begin{bmatrix}
    x\\
    y
\end{bmatrix}=\begin{bmatrix}
    -b_1\\
    -b_2
\end{bmatrix},
\end{equation*}
which, by the Rouch\'e–Capelli Theorem, admits one solution iff $|a|\neq 1$. The case $b=0$ follows directly.

If $|a|=1$, using polar coordinates $z=\rho e^{i\theta}$ and $a=e^{i\varphi}$, we write
\begin{align*}
z+a\bar{z}+b=\rho(e^{i\theta}+e^{i(\varphi-\theta)})+b&=\rho (e^{i(\theta-\varphi/2)}+e^{i(\varphi/2-\theta)})e^{i\varphi/2}+b\\
&=\rho f(\theta)\sqrt{a}+b,
\end{align*}
where $f(\theta)=e^{i(\theta-\varphi/2)}+e^{i(\varphi/2-\theta)}$ is a real valued function with image the set $[-2,2]$. Thus, the map $\rho f(\theta)$ has image $\R$ and the equation $z+a\bar{z}+b=0$ admits an infinite amount of solutions if $b/\sqrt{a}\in\R$ and no solutions otherwise.
\end{proof}

\begin{prop}\label{n=m}
The harmonic equation $z^n+a\bar{z}^n+b=0$ for $a,b\in\C$ admits $n$ solutions iff $|a|\neq 1$ and $b\neq 0$.
If $|a|\neq 1$ and $b=0$ the only solution is $z=0$.

If $|a|=1$ the equation admits an infinite amount of solutions if $b/\sqrt{a}\in\R$ and zero solutions otherwise.
\end{prop}
\begin{proof}
Note that, a solution of the equation $z+a\bar{z}+b=0$ generates for the equation $z^n+a\bar{z}^n+b=0$ the solutions given by the $n$-th roots. 

On the other hand, a solution of the equation $z^n+a\bar{z}^n+b=0$ generates for the equation $z+a\bar{z}+b=0$ the solutions given by the $n$-th power.

The thesis follows from this argument and Lemma~\ref{deg1}.
\end{proof}

Related to this problem we make the following

\begin{con}
Some computed examples suggest that the characteristic polynomial of the matrix representing the Hermitian Killing form with respect to the basis $\lbrace [z^{\alpha}w^{\beta}]\mid\alpha,\beta< n\rbrace$ of the quotient $\C[z,w]/\langle h,h^{\ast}\rangle$ for $h(z,w)\equiv h(z,\bar{z})=z^n+a\bar{z}^n+b$ with $|a|\neq 1$ is
\begin{equation*}
\scriptstyle
  q_{\C}(\lambda)=\frac{(\lambda-n^2)((|a|^2-1)\lambda\pm n^2|b-a\bar{b}|)^{n-1}((|a|^2-1)^2\lambda-n^2|b-a\bar{b}|^2)^{\frac{n(n-1)}{2}}((|a|^2-1)^2\lambda+n^2|b-a\bar{b}|^2)^{\frac{(n-1)(n-2)}{2}}}{(|a|^2-1)^{2n(n-1)}}.
\end{equation*}

If $b\neq 0$, the difference of the number of positive and the number of negative roots of $ q_{\C}(\lambda)$ is
\begin{equation*}
1+(n-1)+\frac{n(n-1)}{2}-\left((n-1)+\frac{(n-1)(n-2)}{2}\right)=n,
\end{equation*}
in accordance to Proposition~\ref{n=m}.

If $b=0$, the polynomial simplifies $q_{\C}(\lambda)=\lambda^{n^2-1}(\lambda-n^2)$ and the difference of the number of positive and the number of negative roots is $1-0=1$, in accordance to Proposition~\ref{n=m}.
\end{con}

%-------------------------------------------------------------------------------------------

\section{Conclusion}
We presented a new object, which could be useful to determine solutions of polynomials with conjugate variables. This may open the way for computations that will lead to new results, such as the one of Theorem~\ref{main}.

%-------------------------------------------------------------------------------------------

\section{Acknowledgments}
The author wants to profoundly thank Prof.\ Giorgio Ottaviani for all the valuable suggestions he gave and the time he spent reading many versions of this work.

%-------------------------------------------------------------------------------------------




\end{document}